\documentclass[12pt,a4paper]{article}
\usepackage[dvips]{graphicx}
\usepackage{setspace}
\usepackage{amsmath}
\usepackage{overcite}
\usepackage{amsbsy}
\usepackage{amssymb}
\usepackage[utf8]{inputenc}
\usepackage[T1]{fontenc}

\newcommand{\ud}{\ensuremath{\mathrm{d}}}

\newcommand{\bmi}{\ensuremath{\boldsymbol}}

\newcommand{\bnabla}{\ensuremath{\boldsymbol\nabla}}
\newcommand{\bcdot}{\ensuremath{\boldsymbol\cdot}}

\newcommand{\p}{\ensuremath{\partial}}
\newcommand{\beq}{\begin{equation}}
\newcommand{\eeq}{\end{equation}}
\newcommand{\eps}{\ensuremath{\epsilon}}

\long\def\symbolfootnote[#1]#2{\begingroup%
\def\thefootnote{\fnsymbol{footnote}}\footnote[#1]{#2}\endgroup}

\begin{document}

\begin{center}
{\bf \Large Phototactic bioconvection in a forward scattering suspension illuminated by both diffuse and oblique collimated flux}\\[12pt]
{M. K. Panda$^1$\symbolfootnote[1]{Corresponding author;
e-mail:mkpanda@iiitdmj.ac.in}, S. K. Rajput$^2$ }\\
{$^1$ Department of Mathematics, PDPM Indian Institute of Information Technology Design and Manufacturing, Jabalpur 482005, India}
\end{center}

\noindent
{\bf Abstract}

The onset of light-induced bioconvection via linear stability theory is investigated qualitatively for a suspension of phototactic algae.
The forward scattering algal suspension is uniformly illuminated by both diffuse and oblique collimated flux. An unstable mode of disturbance at bioconvective instability transits from the stationary (overstable) to overstable (stationary) state at the variation in forward scattering coefficient for fixed parameters. Suspension becomes more stable as forward scattering coefficient increases.

\newpage
 \section{INTRODUCTION}
\label{sec1}

In 1961, Platt coined the field of bioconvection, an important phenomenon that results a self-organized structure and flow in aquatic environments and is usually recognized by patterns in non-neutrally buoyant cell concentration \cite{plat,pk:kp,hp:ph,bees20,javadi20}. The overall picture is that the term bioconvection describes hydrodynamic instabilities, which evolve as a resultant of swimming behaviors and physical properties of the cells, such as density and fluid flows etc. In pattern formation, the dense collections of biased swimming microorganisms (e.g., single-celled algae, bacteria and protozoa etc.) are heavier than surrounding medium (water here) but can swim on average upward in many situations. If the cells are non-motile, then the bioconvection patterns disappear. However, there exist cases where the physical properties of the cells such as up-swimming and higher density are not involved in the process of pattern formation \cite{pk:kp}. For a range of species, the typical patterns in bioconvecion are obtained in the laboratory. However, these have also been found in micropatches of zooplankton \emph{in situ} \cite{kils}. Microorganisms possess the ability and structures that would allow them to propel themselves through their environment and these active or passive motile responses are called \textit{taxes}. For instance, the directed response of cells to a light gradient in their local environment is called  {\it phototaxis}, their swim against gravity is named as {\it gravitaxis} and a combination of viscous and gravitational torques for the bottom-heavy microorganisms can lead them to {\it gyrotaxis}. Phototaxis is relevant to this article only.   

Investigation on pattern formation via bioconvection in an algal suspension reveals that a combination of diffuse and oblique collimaed flux can influence the fluid flows and  the corresponding cell concentration patterns via their shape, size and/or scale etc.\cite{wa,ka85,vin95,williams_11,kh_1997,kage_13,ka86}. In well-stirred collections of photosynthetic microorganisms, the evolved dynamic patterns may be damaged or unbroken by intense light. The steady patterns can be destroyed by strong light in suspensions of microorganisms  or the formation of patterns can be prevented by strong light in well-stirred cultures. Key to the changes in bioconvection patterns via light intensity is algae swim toward (or away from) moderate (or intense) light to photosynthesize (or to avoid damage) \cite{had87}. In addition, absorption and scattering to incident light by the microorganisms may be second reason \cite {mkp20}. Furthermore, the light scattering can be categorized into isotropic and anisotropic scattering depending on the size, shape and the relative refractive index of the algal cells. Isotropic scattering represents  distribution of light uniformly across all directions, whereas the opposite is true in case of anisotropic scattering. Again, anisotropic scattering can be divided into forward and backward scattering. Backward (or forward) scattering scatters more energy into the backward (or forward) directions. In the visible wavelength range, the algae scatter mostly in the forward direction via their size parameter \cite{privoznik,pilon}.

\begin{figure}[!h]
\begin{center}
\includegraphics[scale=0.5]{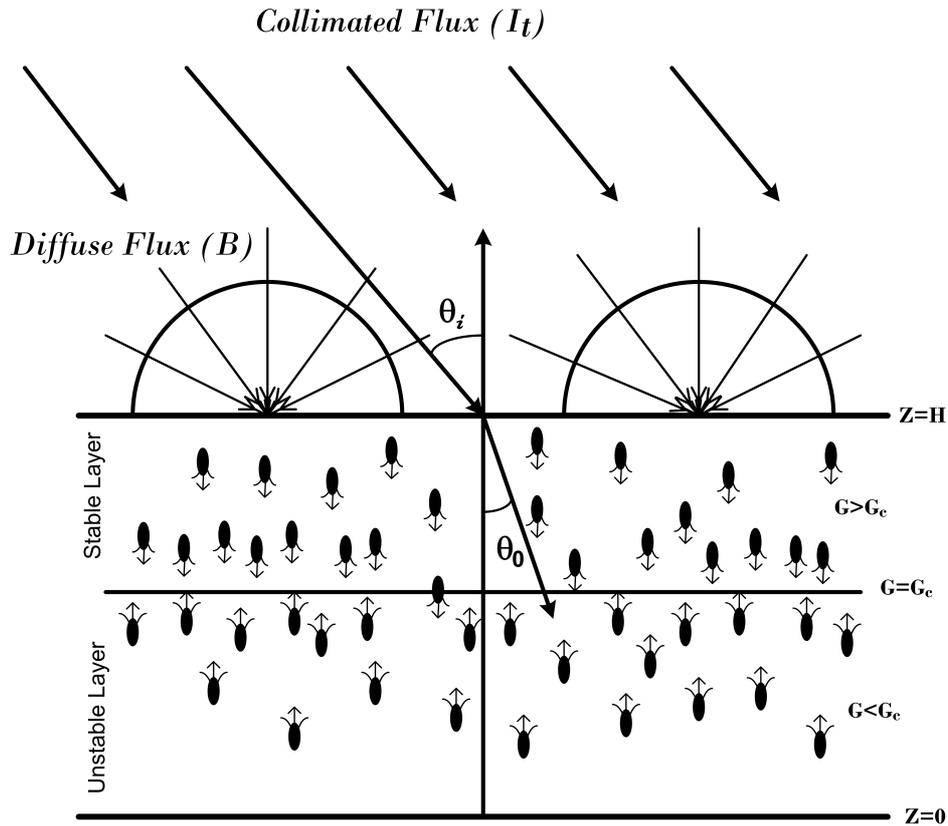}
\end{center}
\caption{$G_c$ is the critical value of total intensity, where algae aggregate.}
\label{fig1}
\end{figure}  

To study this bioconvective system, we employ the phototaxis model proposed by Panda \cite{mkp20}, where the forward scattering algal suspension is illuminated by both diffuse and vertical (normal) collimated flux.  In a natural environment, the sunlight strikes the algal suspension at some off-normal angle via an oblique collimated flux.  Also, the sunlight is able to reach in deeper waters via forward scattering by algae and it may influence the radiation field i.e., light intensity and radiative heat flux etc. Therefore, the light intensity profiles may be redistributed and the swimming behavior (or phototaxis) of algae may be affected if the suspension is illuminated by both diffuse and oblique collimated flux. Looking at the bioconvective instabilities associated with this system, the realistic and reliable models on phototaxis should include the illuminating source to be a combination of diffuse and oblique collimated flux to describe the swimming behavior of algae accurately \cite{bees20,williams_11,incropera,ishikawa}. Furthermore, the generated bioconvection due to both diffuse and oblique collimated flux on an algal suspension should not be disregarded when applying it in industrial exploitation \cite{bees20,javadi20,kils}. Panda  \cite{mkp20} ignored the effects of oblique collimated flux  in addition to the diffuse flux as an illuminating source to the algal suspension and qualitatively analysed the radiation field. In this paper, we include the effect of oblique collimated flux as an illuminating source to algal suspension to investigate the  bioconvective instability in contrast to Panda \cite{mkp20}.

\begin{figure}[!h]
\begin{center}
\includegraphics[scale=0.5]{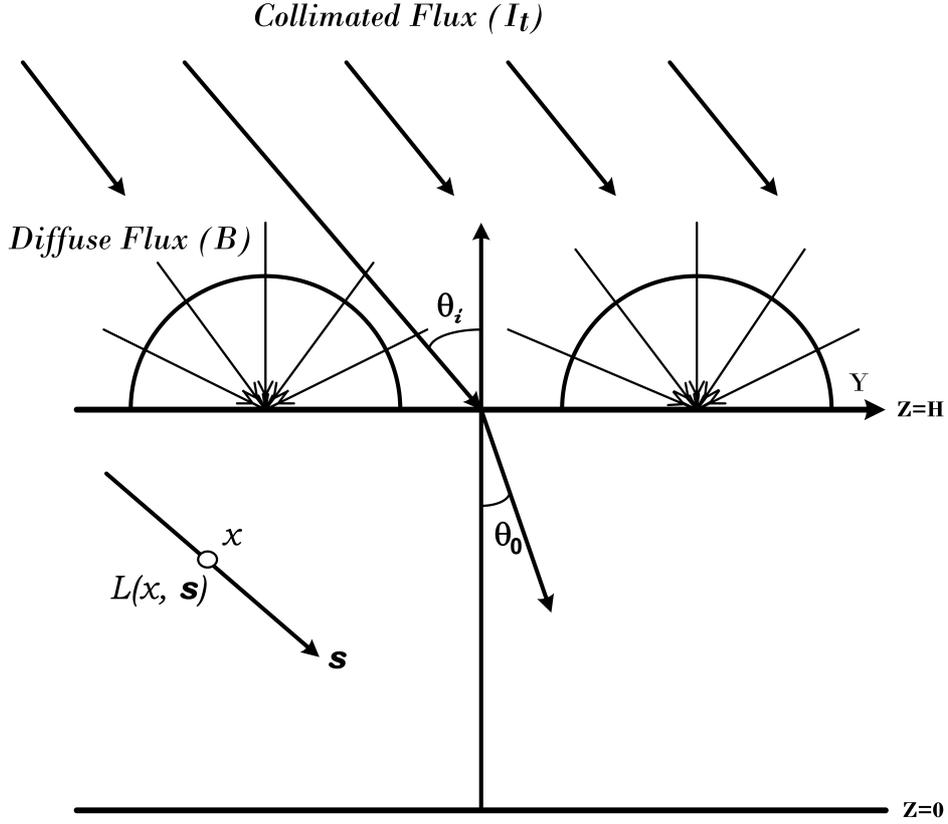}
\end{center}
\caption{Diffuse and oblique collimated flux as an illuminating source to the algal suspension.}
\label{fig2}
\end{figure}  

Look at the bioconvective instabilities associated with these light-induced patterns in a dilute algal suspension. The basic state for zero flow of the bioconvective system is the equilibrium state between phototaxis and cell diffusion. As a result, the cells accumulate as a horizontal concentrated sublayer across the suspension and the sublayer position is a function of critical intensity $G_c$. The sublayer divides the whole domain into two regions based on $G_c$  i.e. a region above the sublayer ($G>G_c$) and a region below it ($G<G_c$) [see Fig. \ref{fig1}]. The fluid motions in the  unstable lower region can penetrate the stable upper region if the system becomes unstable and this mechanism is an example of the penetrative convection \cite{bst}.

There has been significant progress in modeling phototaxis that result in bioconvection and various approaches to them, through linear stability theory and numerical simulation have been published. Vincent and Hill \cite{vh:hv} investigated the onset phototactic bioconvection and found non-oscillatory and non-stationary modes of disturbance. Two-dimensional bioconvective flow pattern was simulated by Ghorai and Hill \cite{ghp} via the model of Vincent and Hill \cite{vh:hv}. But both of these studies neglected scattering by algae. Ghorai \textit{et al.} \cite{gph} investigated the onset of bioconvection in an isotropic scattering suspension via a linear stability theory. However, a steady-state profile with sublayer located at two different depths has been noticed in their study due to isotropic scattering by algae for some parameters. Ghorai and Panda \cite{gp} examined linear stability analysis in a forward scattering algal suspension. At bioconvective instability, they reported a transition from a non-oscillatory  (non-stationary) to an non-stationary (non-oscillatory) mode due to forward scattering. Panda and Ghorai \cite{pg} numerically computed bioconvection for an isotropic scattering suspension confined in a two-dimensional geometry. The evolved bioconvective flow found by them was different to Ghorai and Hill \cite{ghp} qualitatively. Afterwards, Panda and Singh \cite{pr} simulated two-dimensional phototactic bioconvection in between rigid sidewalls via the model of Vincent and Hill \cite{vh:hv}. The lateral rigid walls enhance the stabilizing effect on suspension as reported by them. However, these studies did not include diffuse flux as an illuminating source. Panda \textit{et al.}  \cite{prmm16}  investigated the impact of vertical collimated flux and diffuse flux together on phototactic bioconvection and concluded that diffuse flux makes the suspension more stable. Furthermore, the bimodal vertical base concentration profiles switch to unimodal ones at the effects of diffuse flux. Panda \cite{mkp20} investigated the forward scattering effect on the bioconvective instability with diffuse and vertical collimated flux together. He reported that the bimodal vertical base concentration profiles transit into the unimodal ones and vice versa due to forward scattering when scattering is strong. However, the effects of oblique collimated flux were not incorporated in the aforementioned reports. First time, Panda \textit{et al.} \cite{PPS} incorporated the oblique collimated flux to investigate their effect on bioconvection in a non-scattering algal suspension. They predicted that the bioconvective solutions switch from a non-oscillatory state to overstable state and vice versa at bioconvective instability due to oblique collimated flux.  Furthermore, the mode $2$ instability transits to a mode $1$ instability due to oblique collimated flux as reported by them. More recently, Kumar \cite{sandeep} investigated the impact of oblique collimated flux on bioconvective solutions in an isotropic scattering algal suspension. He found that the bioconvection solutions become mostly oscillatory  in contrast to stationary in the case of a rigid upper boundary. However, phototactic bioconvection with oblique collimated and diffuse flux together on a forward scattering algal suspension is still not fully explored. Therefore bioconvective instability on a forward scattering algal suspension using a realistic phototaxis model is investigated here. 

The structure of the paper is as follows. We start with the mathematical formulation of the proposed work. Next, an equilibrium solution of the full equations is found and is then perturbed. The perturbed problem is solved next by using linear stability theory. Finally the numerical results are summarized and novelty of the proposed model is discussed.

\section{Mathematical formulation}
\label{sec2}

We consider a forward scattering algal suspension, filling the space between two parallel boundaries of infinite extent in the $(y,z)$ plane [see Fig. \ref{fig2}]. Reflection from the top and the bottom boundaries is ignored and both the oblique collimated and diffuse flux illuminate uniformly the upper surface of the forward scattering algal suspension. The oblique collimated flux illuminates the forward scattering algal suspension and then, it is transmitted across the suspension  [see Fig. \ref{fig2}].  Since the refractive index of the water is different from that of air, the refraction angle $\theta_{0}$ is related to the incidence angle $\theta_{i}$ by Snell's law, $\sin{\theta_{i}}=n_{0}\,\sin{\theta_{0}},$ where $n_{0}$ is the index of refraction of the water. The estimated value for refractive index of the water is $1.333$ approximately \cite{Daniel_1979} and the index of refraction of air has been assumed to be equal to unity.  The light incident on the suspension is absorbed and scattered by algae in forward direction due to difference in the index of refraction between algae and water.

\subsection{Governing equations}
\label{subsec2_1}

Based on the previous continuum models \cite{pk:kp,mkp20}, each cell has a volume $\vartheta$ and density 
$\rho+\delta\rho$, where $\rho$ is the density of the fluid in which the 
cells swim and $\delta\rho/\rho \ll 1$. Let ${\bmi u}$ and $n$, respectively,
denote the velocity and concentration in the suspension. Accounting for suspension incompressibility condition, the conservation of mass gives  

\begin{equation}
\bnabla\bcdot{\bmi u} = 0.  \label{eq4}
\end{equation}

The momentum equation via Boussinesq approximation is given by \cite{mkp20}

\begin{equation}
\rho\frac{\mathrm{D}{\bmi u}}{\mathrm{D} t} = -\bnabla P_e + \mu{\nabla^2}{\bmi u} -
n\,\delta\rho\, g \vartheta \hat {\bmi z}, \label{eq5}
\end{equation}

where $P_e,$ $g,$ and $\mu$ are the excess pressure, acceleration due to gravity and viscosity of the suspension (equal to fluid) respectively. Also, $\mathrm{D}/\mathrm{D}t=\p/\p t + {\bmi u}\bcdot\bnabla$ denotes the convective time derivative.

The governing equation for cell conservation is

\begin{equation}
\frac{\p n}{\p t} = -{\bmi \nabla}\bcdot{\bmi F}, \label{eq6}
\end{equation}

where  cell flux is given by

\begin{equation}
{\bmi F} = n {\bmi u}+{n} \mathbf{W}_c - \mathsf{\mathbf D} {\bmi \nabla}{n}.  \label{eq7}
\end{equation}
We choose $\mathsf{\mathbf D}=D \mathsf{\mathbf I}$[see Panda \cite{mkp20}]. Thus, 

\beq
{\bmi F} = n {\bmi u}+{n} \mathbf{W}_c - D {\bmi \nabla}{n}. \label{xflux}
\eeq
Two key assumptions are considered in this study as similar to Panda \cite{mkp20} so that the Fokker-Planck equation can be eliminated from the full equations of phototactic bioconvective system.

\subsection{The mean swimming orientation}
\label{subsec2_2} 

The governing equation to express light intensity i.e. the  radiative transfer equation (hereafter referred to as RTE) is given by \cite{modest,chandras}

\beq
{\bmi s}\bcdot\bnabla L({\bmi x},{\bmi s}) + (\kappa+\sigma)L({\bmi x},{\bmi s}) =
\frac{\sigma}{4\pi}\int_0^{4\pi}L({\bmi x},{\bmi s}') p({\bmi
s}',{\bmi
s})\,\ud\Omega', \label{eq1}
\eeq   

where $\Omega',$ $\kappa,$ and $\sigma,$  are the  solid angle, absorption coefficient, the scattering coefficient  respectively. The scattering phase function $p({\bmi s}',{\bmi s}),$ is assumed to be linearly-anisotropic \cite{modest}. Thus,

\beq
p({\bmi s}',{\bmi s})=1+\mathrm{A_1}\cos{\theta}\cos{\theta'},  \label{phse}
\eeq

where $\mathrm{A_1}$ is the linearly-anisotropic scattering coefficient, which indicates backward ($-1 < \mathrm{A_1} <0$), forward ($0<\mathrm{A_1} \le 1$) or isotropic scattering ($\mathrm{A_1} =0$) respectively.

The light intensity at a top location i.e. ${\bmi x}_b=(x,y,H)$ of the algal suspension is given by

$$
L({\bmi x}_b,{\bmi s}) = \mathrm{I}_t \,\delta({\bmi s}-{\bmi s}_0)+\dfrac{B}{\pi} = \mathrm{I}_t\,
\delta(\cos\theta-\cos\theta_0)\delta(\phi-\phi_0)+\dfrac{B}{\pi},
$$

where $\mathrm{I}_t, $ and $B$ are the magnitudes of oblique collimated and diffuse flux respectively [refer Panda \cite{mkp20} for details ]. Also, the incident direction is given by ${\bmi s}_0=\cos(\pi-\theta_0) \hat{\bmi z} + \sin(\pi-\theta_0) (\cos\phi_0  \hat{\bmi x} + \sin\phi_0 \hat{\bmi y}),$  where $\hat {\bmi x}, \hat{\bmi y}, \hat {\bmi z}$ are unit vectors along the $x,y,z$ axes. 

We assume $\kappa = \varkappa n$ and $\sigma=\varsigma n$ and in terms  of the scattering albedo $\omega=\sigma/(\kappa+\sigma),$ the RTE  in a forward scattering algal suspension thus becomes
\beq
{\bmi s}\bcdot\bnabla L({\bmi x},{\bmi s}) + \beta  L({\bmi
x},{\bmi s}) =\frac{\omega\beta }{4\pi}\int_0^{4\pi}L({\bmi
x},{\bmi s}')(1+\mathrm{A_1}\cos{\theta}\cos{\theta'})\,\ud\Omega'. \label{eq2}
\eeq
Here the extinction coefficient is given by $\beta=(\varkappa+\varsigma) n$. Also, $0 \le \omega \le 1$ and the medium is purely absorbing (or scattering) when $\omega=0$ (or $\omega=1$).

The total intensity, $G({\bmi x}),$ at a point ${\bmi x}$ in the medium is 

$$
G({\bmi x})=\int_0^{4\pi}{L({\bmi x},{\bmi s})}\,\ud\Omega,
$$ 

and the radiative heat flux, ${\bmi q}({\bmi x}),$ is 

$$
{\bmi q}({\bmi x}) = \int_0^{4\pi}{L({\bmi x},{\bmi s})}\,{\bmi s}\,\ud\Omega.
$$  

The average swimming direction, $<{\bmi p}>$ is defined as the ensemble-average of the swimming direction, ${\bmi p}$ for all the cells in a small volume.   The average swimming velocity, $\mathbf{W}_c,$ is thus

$$
\mathbf{W}_c = {W_c}<{\bmi p}>,
$$
where the ensemble average swimming speed of cells is given by $W_c$.

The mean swimming direction, $<{\bmi p}>$, is given by \cite{mkp20}

\beq
<{\bmi p}> = -M(G) \frac{{\bmi q}}{|{\bmi q}|}, \label{eq3}
\eeq

where  $M(G)$ is the phototaxis function
such that

$$
M(G)\quad
\left\{ \begin{array}{c}
\ge 0 \qquad \mbox{if } G \le G_c, \\
< 0 \qquad \mbox{if } G > G_c.\\
\end{array}
\right.     
$$

\subsection{Boundary conditions}
\label{subsec2_3}

The bottom (top) boundary is taken as rigid i.e. $z=0$ (free i.e. $z=H$) and zero cell flux through them. Thus, 

\beq
{\bmi u}= {\bmi F}\bcdot {\hat {\bmi z}} = 0\qquad \mbox{on}\quad z=0,\label{eq8}\\
\eeq
\beq
\frac{\p^2}{\p z^2}({\bmi u}\bcdot\hat {\bmi z})= {\bmi u}\bcdot \hat {\bmi z}=
{\bmi F}\bcdot {\hat {\bmi z}}=0\qquad\mbox{on}\quad z=H.\label{eq9}
\eeq

 The nonreflecting top and the bottom  boundary conditions are given by

\begin{eqnarray}
\mbox{at }z=H,\, L(x,y,z,\theta,\phi) &=& \mathrm{I}_t
\,\delta(\cos\theta-\cos\theta_0)\delta(\phi-\phi_0)+\dfrac{B}{\pi},\, \pi/2\le\theta\le
\pi, \label{eq10}\\
\mbox{at }z=0,\, L(x,y,z,\theta,\phi) &=& 0,\, \qquad 0\le\theta\le \pi/2.
\label{eq11}
\end{eqnarray}                               

\subsection{Scaling of the equations}
\label{subsec2_4}  

 After scaling the length by $H,$ time by $H^2/D,$ the bulk fluid velocity by $D/H,$ the pressure by $\mu D/H^2,$ the cell
concentration by $\bar n,$ and intensity by $\mathrm{I}_t,$ the governing bioconvective system becomes

\begin{eqnarray}
&& \bnabla\bcdot{\bmi u} = 0, \label{eq12}\\
&& S_c^{-1}\frac{\mathrm{D}{\bmi u}}{\mathrm{D} t} = -\bnabla P_e + {\nabla^2}{\bmi u} -
R\,n \hat {\bmi z}, \label{eq13} \\
&& \frac{\p n}{\p t} = -{\bmi\nabla}\bcdot{\bmi F}, \label{eq14}
\end{eqnarray}

where the cell flux  is given by

$$
{\bmi F} = n {\bmi u}+{n}V_c <{\bmi p}> - {\bmi\nabla}{n}.  
$$

Here $S_c=\mu/\rho D$ is the Schmidt number, 
$V_c=W_c H/D$ is the scaled swimming speed, and $R=\bar n \vartheta  g
H^3\delta\rho/\mu D$ is the Rayleigh number. 

In terms of the nondimensional variables, the 
RTE (see Eq. (\ref{eq2})) becomes 

\beq
{\bmi s}\bcdot\bnabla L({\bmi x},{\bmi s}) + \tau_H n L({\bmi x},{\bmi s}) =
\frac{\omega\tau_H n}{4\pi}\int_0^{4\pi}L({\bmi x},{\bmi s}')(1+\mathrm{A_1}\cos\theta\,\cos\theta') \,\ud\Omega',  \label{aret}
\eeq

where $\tau_H=(\varkappa+\varsigma) \bar n H$ is the (vertical) optical depth of
the suspension. In terms  of the direction cosines $(\xi,\eta,\nu)$ of the 
unit vector ${\bmi s}$, where

$$
\xi=\sin\theta\cos\phi,\quad \eta=\sin\theta\sin\phi,\quad \nu=\cos\theta,
$$  

the nondimensional RTE (see Eq. (\ref{aret})) can be written as \cite{modest} 

\beq
\xi\frac{\p L}{\p x} + \eta \frac{\p L}{\p y} +\nu \frac{\p L}{\p z} + \tau_H n
L({\bmi x},{\bmi s}) =
\frac{\omega\tau_H n}{4\pi}\int_0^{4\pi}L({\bmi x},{\bmi s}')(1+\mathrm{A_1}\cos\theta\,\cos\theta') \,\ud\Omega'.
\label{eq15}
\eeq      

Please note that we have used the same notation for the dimensional and nondimensional variables here. A typical photoresponse curve i.e. $M(G)$\cite{mkp20} can be expressed as:                    

\begin{equation}
M(G)=0.8\sin{\left(\frac{3\pi}{2}\chi{\left(G\right)}\right)}-0.1\sin{\left(\frac{\pi}{2}\chi{\left(G\right)}\right)},\quad
\chi{\left(G\right)}=\dfrac{G}{3.8} e^{\Upsilon\left(3.8-G\right)},   
\label{eq0}
\end{equation}
where $\Upsilon=0.252$  and the critical intensity  occurs at location $G=G_c \approx 1.3$ [see Fig. \ref{fig4}(a)].

\section{The basic solution}
\label{sec3}

An equilibrium solution is found via setting

\beq
{\bmi u}=0,\,n=n_s(z),\,P_e=P_e{^s},\quad\mbox{and}\quad L=L_s(z,\theta), \label{eq16}
\eeq

in Eqs.(\ref{eq12})--(\ref{eq14}) and Eq. (\ref{eq15}).

At the basic state, the total intensity $G_s(z)$ and the radiative heat flux ${\bmi q}_s(z)$  are given by

$$
G_s(z)=\int_0^{4\pi}L_s(z,\theta)\,\ud\Omega,\quad {\bmi q}_s(z) =
\int_0^{4\pi}L_s(z,\theta)\,{\bmi s}\,\ud\Omega.
$$    
     
The $x$ and $y$ components of ${\bmi q}_s$ vanish since $L_s(z,\theta)$ does not depend on $\phi$. Therefore, ${\bmi q}_s=-q_s\hat{\bmi z}$, where $q_s=|{\bmi q}_s|$. From Eq. (\ref{eq15}), the basic state intensity $L_s$ follows as:

\beq
\frac{\p L_s}{\p z} + \frac{\tau_H n_s L_s}{\nu} =\frac{\omega\tau_H
n_s}{4\pi\nu}{\Big(} G_s(z)-\mathrm{A_1} q_{s} \nu{\Big)}. \label{eq17} 
\eeq 
 
We take $L_s = L_s^c + L_s^d$. Then  $L_s^c$ satisfies

\beq
\frac{\p L_s^c}{\p z} + \frac{\tau_H n_s L_s^c}{\nu} =0, \label{col1}                
\eeq    
 
subject to the top boundary condition

\beq
L_s^c(1,\theta) =  \delta(\nu-\nu_0)\delta(\phi-\phi_0). \label{col2}
\eeq

Solving Eqs. (\ref{col1}) and (\ref{col2}), we get 
          
$$
L_s^c=\exp\left(\int_z^1\frac{\tau_H
n_s(z')}{\nu}\,\ud z'\right)\delta(\nu-\nu_0)\delta(\phi-\phi_0).
$$

Now the diffused component $L_s^d$ satisfies

\beq
\frac{\p L_s^d}{\p z} + \frac{\tau_H n_s L_s^d}{\nu} =\frac{\omega\tau_H
n_s}{4\pi\nu}{\Big(}G_s(z)-\mathrm{A_1} q_s \,\nu{\Big)},  \label{eq18}               
\eeq      
 
subject to the boundary conditions

\begin{eqnarray}
L_s^d(1,\theta) &=& \dfrac{B}{\pi},\quad \pi/2\le\theta\le \pi,\label{bcin1}\\
L_s^d(0,\theta) &=& 0,\quad  0\le\theta\le \pi/2. \nonumber
\end{eqnarray}

We define a new variable 

$$
\tau = \int_z^1 \tau_H n_s(z')\,\ud z',
$$ 

so that the total intensity $G_s$ and radiative heat flux ${\bmi q}_s$ become 
function of $\tau$ only.

The basic total intensity is written as

\begin{equation}
G_s=G_s^c+G_s^d,
\label{eq22}
\end{equation}

 where

$$
G_s^c=\int_0^{4\pi}L_s^c(z,\theta)\,\ud\Omega=\exp\left(-\tau_H/\cos{\theta_{0}}\int_z^1
n_s(z')\,\ud z'\right),$$ and
$$
G_s^d =\int_0^{4\pi}L_s^d(z,\theta)\,\ud\Omega.  
$$

Similarly, 

\beq
{\bmi q}_s =  {\bmi q}_s^{c}+{\bmi q}_s^{d}, \label{eq23} 
\eeq

where

$$ 
{\bmi q}_s^{c} = \int_0^{4\pi}L_s^c\,{\bmi s}\,\ud\Omega=-\cos{\theta_{0}} \exp
\left(-\tau_H/\cos{\theta_{0}}\int_z^1 n_s(z')\,\ud z'\right)\hat{\bmi z},$$ and
$$
{\bmi q}_s^{d} = \int_0^{4\pi}L_s^d(z,\theta)\,{\bmi s}\,\ud\Omega. 
$$
Eqs. (\ref{eq22}) and (\ref{eq23}) recast as two set of coupled Fredholm integral equations of the second
kind \cite{modest} via variable $\tau$ i.e.

\begin{eqnarray}
&&G_s(\tau) = 2 B E_{2}(\tau)+e^{-\tau / \cos{\theta_{0}}} \nonumber\\&&+
\frac{\omega}{2}\int_0^{\tau_H}\Big(
G_s(\tau')\,E_1\left(|\tau-\tau'|\right) +
\mathrm{A_1} \,\textnormal{sgn}(\tau-\tau')q_s(\tau')E_2\left(|\tau-\tau'|\right)\Big)\ud\tau',
\label{eq24} \quad \quad \\[8pt]
&&q_s(\tau)=2 B E_{3}(\tau)+\cos{\theta_{0}} \, e^{-\tau / \cos{\theta_{0}}}\nonumber\\&&+\frac{\omega}{2}
\int_{0}^{\tau_H}\Big(\mathrm{A_1}\,q_s(\tau')E_3\left(|\tau-\tau'|\right)+\textnormal
{sgn}(\tau-\tau')G_s(\tau')E_2(|\tau-\tau'|)\Big)\ud\tau'. \quad  
\label{eq25}
\end{eqnarray}

Here, `sgn' is the sign function and $E_n(x)$ is the exponential integral function of order $n$ \cite{chandras}. Using the method of subtraction of singularity \cite{press}, Eqs. (\ref{eq24}) and (\ref{eq25}) are solved.

\begin{figure}[!h]
\begin{center}
\includegraphics[scale=0.9]{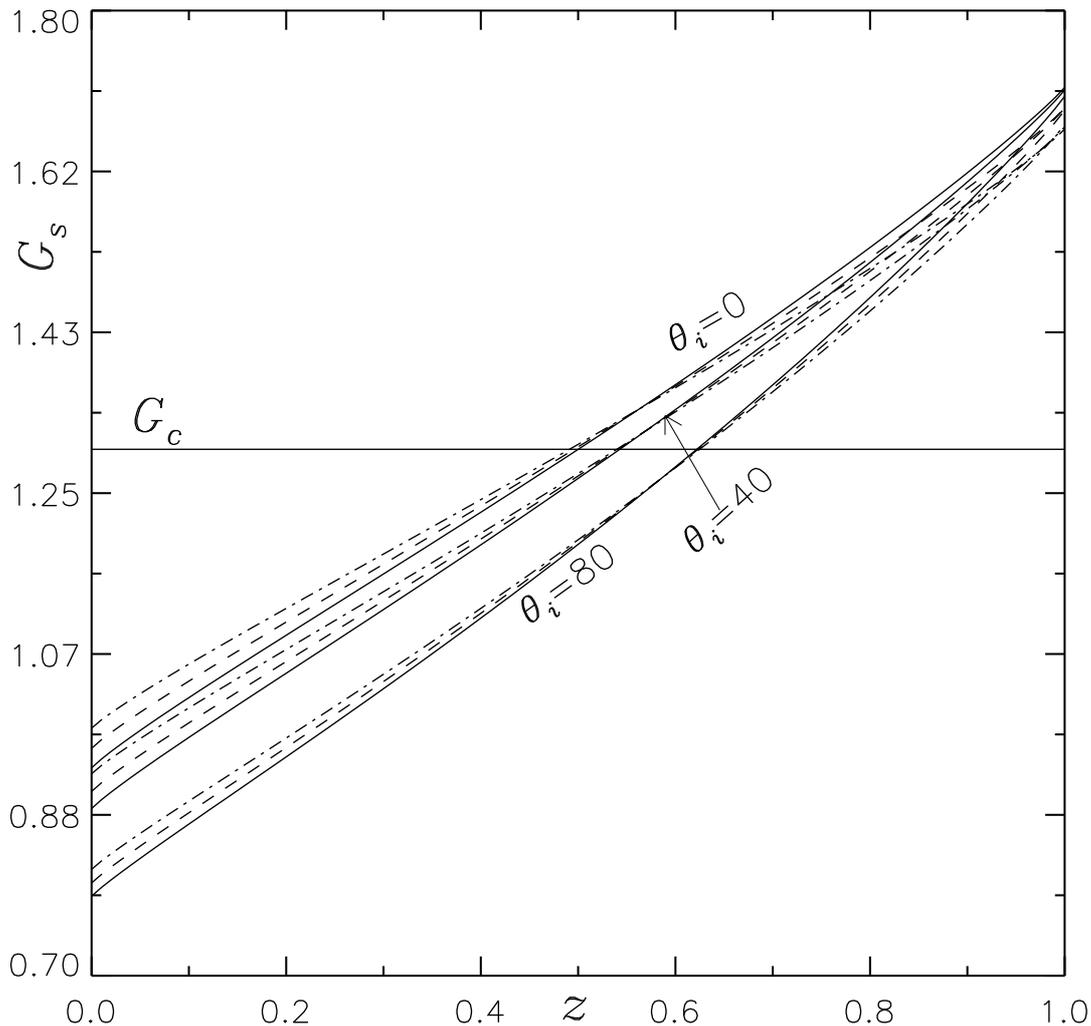}
\end{center}
\caption{Effects of the linearly-anisotropic scattering coefficient $\mathrm{A_1}$ on the total intensity $G_s$ in a uniform suspension for single scattering albedo $\omega=0.4,$ diffuse irradiation $B=0.26$ and different values of angle of incidence $\theta_{i}$  : $\mathrm{A_1}=0$ ($-\-$), $\mathrm{A_1}=0.4$ (-\,-\,-\,-), and $\mathrm{A_1}=0.8$ ($-\!\cdot\!-\!\cdot\!-$). The intersection of the horizontal solid line with $G_s$ represents the location of the critical total intensity $G_c=1.3.$ The fixed parameter is optical depth  $\tau_H=0.5$.}
\label{fig3}
\end{figure}

In the basic state, the average swimming direction becomes 

$$
<{\bmi p}_s>=-M_s\frac{{\bmi q}_s}{q_s}=M_s\hat{\bmi z},
\qquad\quad\mbox{where}\quad
M_s= M(G_s).
$$ 

The concentration, $n_s(z)$, at basic state satisfies 

\beq
\frac{\ud n_s}{\ud z} = V_c\,M_s\,n_s, \label{eq26}
\eeq

subject to

\beq
\int_{0}^{1} n_s\,\ud z = 1.  \label{eq27}
\eeq

A shooting method is used to solve the boundary value problem via Eqs. (\ref{eq24})--(\ref{eq27}) \cite{press}. 

We employ the following set of parameters in accordance with Panda \textit{et al.} \cite{PPS} and Panda \cite{mkp20}:
$S_c = 20,$ $B=0.26,0.5,$ $\omega=0.4,$ $\theta_{i}=0,40,80,$  $\tau_{H}=0.5,1,$ $V_c =15$ and  $\mathrm{A_1}= 0,0.4,0.8$ to explain the equilibrium solution. We fix $S_c = 20,$ $B=0.26,$ $\tau_{H}=0.5,$ and $\omega=0.4$ here and choose the angle of incidence $\theta_{i}$ such that the sublayer at basic state occurs at around mid-height, three-quarter height and top respectively. Next, we vary the forward scattering coefficient $\mathrm{A_1}$ to their effect on total intensity and basic equilibrium solution. Figure \ref{fig3} depicts how the total intensity $G_s$ varies across the depth for a uniform suspension as the forward scattering coefficient $\mathrm{A_1}$ is increased. At the bottom (respectively, top) half of the uniform suspension, the height of total intensity $G_s$ increases (respectively, decreases) as the forward scattering coefficient $\mathrm{A_1}$ is increased. Similarly, sublayer location for a uniform suspension shifts toward the bottom  (respectively, top) across the bottom (respectively, top) half of the domain as the forward scattering coefficient $\mathrm{A_1}$ is increased [see Fig. \ref{fig3}].

\begin{figure}[!h]
\begin{center}
\includegraphics[scale=0.9]{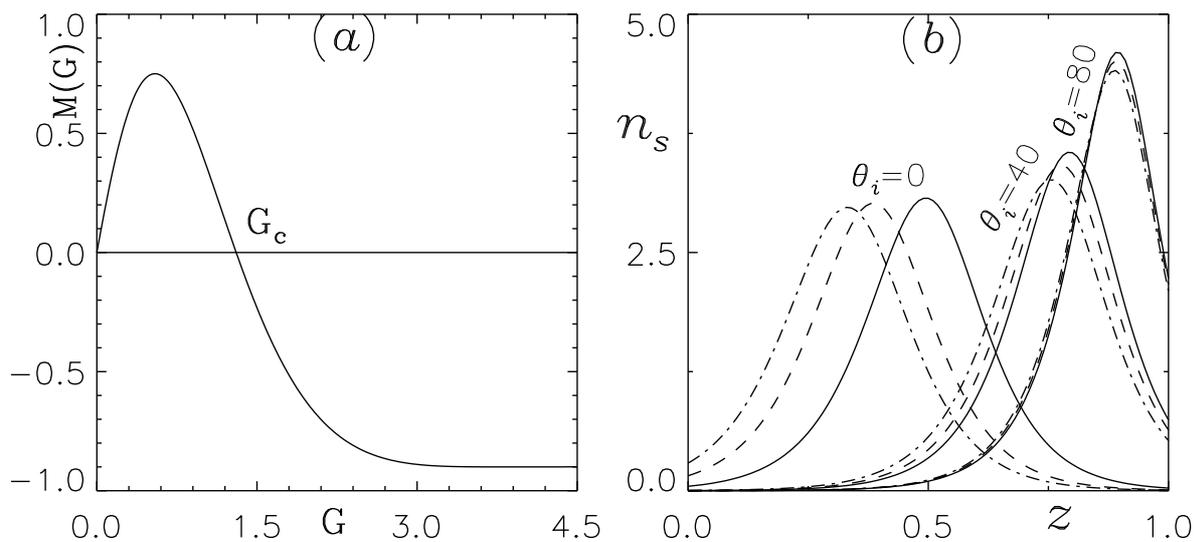}
\end{center}
\caption{(a) The phototaxis function with critical intensity $G_{c}=1.3$ [refer Eq. (\ref{eq0})] and (b) Effects of the linearly-anisotropic scattering coefficient $\mathrm{A_1}$ on the base concentration profiles for  single scattering albedo $\omega=0.4$ and  diffuse irradiation $B=0.5$ : $\mathrm{A_1}=0$ ($-\-$), $\mathrm{A_1}=0.4$ (-\,-\,-\,-), and $\mathrm{A_1}=0.8$ ($-\!\cdot\!-\!\cdot\!-$). Here the fixed parameters are cell swimming speed $V_c=15,$ and optical depth $\tau_H=1$.}
\label{fig4}
\end{figure}

Figure \ref{fig4}(a) shows a photoresponse curve for light intensity with a critical value  $G_{c}=1.3$ and Fig. \ref{fig4}(b) depicts the effect of forward scattering coefficient $\mathrm{A_1}$ on the sublayer at equilibrium state when the parameters $\omega=0.4,$ $B=0.5,$ $V_c=15,$ $\theta_{i}=0,40,80$ and $\tau_{H}=1$ are kept fixed. We start with the case when $\theta_{i}=0.$ Here the sublayer at equilibrium state forms at around mid-height of the suspension domain. In this case, the sublayer location at equilibrium state shifts toward the bottom as $\mathrm{A_1}$ is increased from $0$ to $0.8$ respectively [see Fig. \ref{fig4}]. Similarly, for $\theta_{i}=40$ and $80,$ the sublayer at equilibrium state forms at around three-quarter height and around top of the suspension depth respectively. As $\mathrm{A_1}$ increases from $0$ to $80,$ the sublayer formed at equilibrium state in each case shifts toward the bottom of the suspenion [see Fig. \ref{fig4}].

\begin{figure}[!h]
\begin{center}
\includegraphics[scale=0.95]{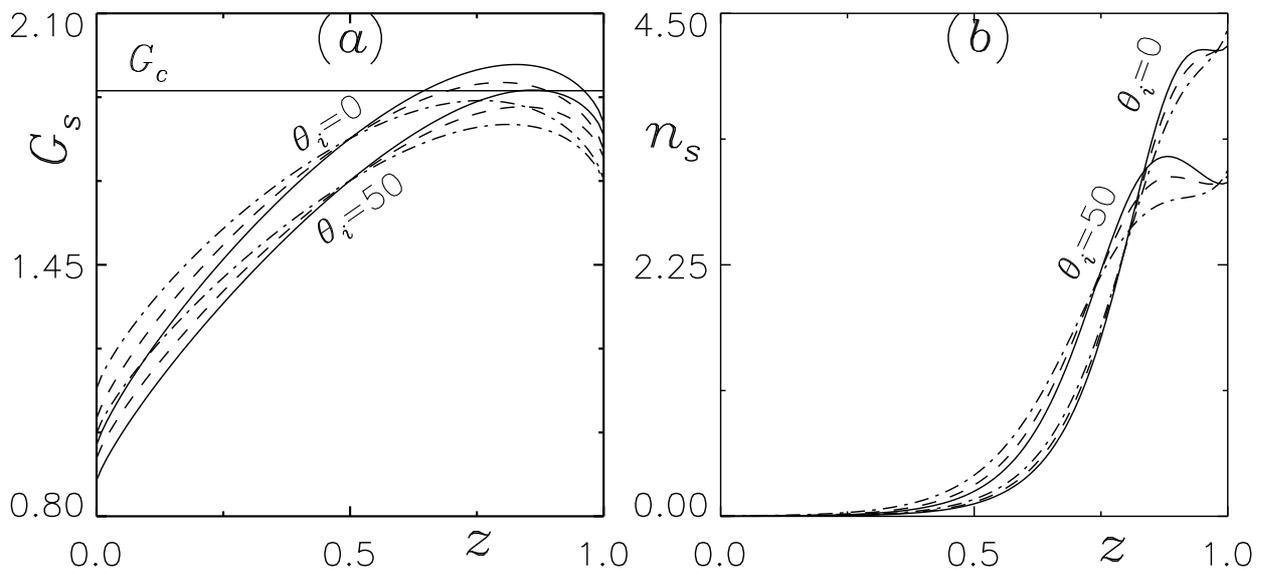}
\end{center}
\caption{(a) Effects of the linearly-anisotropic scattering coefficient $\mathrm{A_1}$ on the total intensity $G_s$ in a uniform suspension for single scattering albedo $\omega=1,$ the cell swimming speed $V_c=15,$ the angle of incidence $\theta_{i}=0,50,$ optical depth $\tau_H=1,$  and diffuse irradiation $B=0.02,$ where the intersection of the horizontal solid line with $G_s$ represents the location of the critical total intensity $G_c=1.9$ [refer Eq. (\ref{eq0}) for $\Upsilon=0.135$] and (b) the corresponding base concentration profiles for cell swimming speed $V_c=15$ : $\mathrm{A_1}=0$ ($-\-$), $\mathrm{A_1}=0.4$ (-\,-\,-\,-), and $\mathrm{A_1}=0.8$ ($-\!\cdot\!-\!\cdot\!-$).}
\label{fig5}
\end{figure}

Figure \ref{fig5}(a) shows the variation of total intensity $G_s$ across the depth of the suspension at the increment of $\mathrm{A_1},$ when the parameters $V_c=15,$ $\theta_{i}=0,50,$ $\tau_H=1,$ $V_c=15,$ $\omega=1,$ $B=0.02$ are kept fixed. The light intensity with a critical value $G_c=1.9$ is selected here. First consider the case when $\theta_{i}=0.$ In this instance, the line for critical intensity intersects at two points to the graph of $G_s$ for $\mathrm{A_1}=0$ and $0.4.$ But, the critical intensity $G_c=1.9$ does not intersect the graph of $G_s$ when $\mathrm{A_1}=80.$ Eventually, the sublayer at equilibrium state forms at two different locations across the suspension for  $\mathrm{A_1}=0$ and $0.4$ and the dual location of the sublayer at equilibrium state shifts to a single location when $\mathrm{A_1}=80$ [see Fig. \ref{fig5}(b)]. 

Next, consider the case when $\theta_{i}=50.$ In this case, the line for critical intensity $G_c=1.9$ intersects at two points to the graph of $G_s$ for $\mathrm{A_1}=0.$  But, the critical intensity $G_c=1.9$ does not intersect the graph of $G_s$ when $\mathrm{A_1}=40, 80$ [see Fig. \ref{fig5}(a)]. Eventually, the sublayer at equilibrium state forms at two different locations across the suspension for  $\mathrm{A_1}=0$ and the dual location of the sublayer at equilibrium state shifts to a single location when $\mathrm{A_1}=40, 80$ [see Fig. \ref{fig5}(b)].

\section{The linear stability problem}
\label{sec4}

The equilibrium solution is perturbed as:

\begin{equation*}
{\bmi u}=0+\epsilon {\bmi u}_1+\mathcal{O}(\epsilon^2),\,n=n_s(z)+\epsilon n_1+\mathcal{O}(\epsilon^2),\,{\bmi
p}={\bmi p}_s+ \epsilon <{\bmi p}_1> +\mathcal{O}(\epsilon^2),\, P_e=P_e{^s}+\epsilon p_e+\mathcal{O}(\epsilon^2), 
\end{equation*}
where ${\bmi u}_1=(u_1,v_1,w_1)$, and $0<\eps\ll 1.$ 

The perturbed governing system is given by

\beq
 \bnabla\bcdot{\bmi u}_1=0,\label{eq28} \\
\eeq

\beq
 S_c^{-1} \frac{\p {\bmi u}_1}{\p t} = -\bnabla p_e - R n_1 \hat{\bmi z} +
\nabla^2{\bmi u}_1, \label{eq29} \\
\eeq

\beq
 \frac{\p n_1}{\p t} + w_1 \frac{\ud n_s}{\ud z} +V_c\,\bnabla \bcdot(<{\bmi
p}_s> n_1 + <{\bmi p}_1> n_s)=\nabla^2 n_1. \label{eq30}
\eeq

Let $L_1$ be the perturbed radiation intensity, i.e., $L=L_s + \eps L_1+\mathcal{O}(\epsilon^2)$. If 
$G=G_s+\eps G_1+\mathcal{O}(\epsilon^2)$, then

\beq
G_1 = \int_0^{4\pi}L_1({\bmi x},{\bmi s})\,\ud\Omega.\label{eq31}
\eeq

Similarly, for the radiative heat flux ${\bmi q}={\bmi q}_s+\eps {\bmi q}_1+\mathcal{O}(\epsilon^2)$, 
we  find

\beq
{\bmi q}_1 = \int_0^{4\pi}L_1({\bmi x},{\bmi s})\,{\bmi s}\,\ud\Omega.\label{eq32}
\eeq

Now the expression

$$
-M(G_s+\eps G_1+\mathcal{O}(\epsilon^2))\frac{{\bmi q}_s + \eps {\bmi q}_1+\mathcal{O}(\epsilon^2)}{|{\bmi q}_s + \eps {\bmi q}_1+\mathcal{O}(\epsilon^2)|}-M_s\hat{\bmi z},
$$

on simplification at $\mathcal{O}(\eps)$, gives perturbed swimming direction

\beq
<{\bmi p}_1>=G_1\frac{\ud M_s}{\ud G} \hat{\bmi z}-M_s\frac{{\bmi q}_1^H}{q_s}, 
\label{eq33}
\eeq

where ${\bmi q}_1^H$ is the horizontal component of the perturbed radiative heat flux
 ${\bmi q}_1$.

Eqs. (\ref{eq28})--(\ref{eq30}) reduce to two equations for $w_1$ and $n_1$ after eliminating $p_e$ (if we take curl of Eq. (\ref{eq29}) twice) and the horizontal components of ${\bmi u}_1$. Accounting the normal mode analysis, we can write
$$
w_1=W(z) \exp(\gamma t + i(lx+my)),\quad n_1=\Theta(z)\exp(\gamma t + i(lx+my)). 
$$
From Eq. (\ref{eq15}), the perturbed intensity $L_1$ satisfies

\beq
\xi \frac{\p L_1}{\p x} + \eta \frac{\p L_1}{\p y} + \nu \frac{\p L_1}{\p
z} + \tau_H n_s L_1 =\frac{\omega\tau_H}{4\pi}\Big(n_s G_1 + G_s
n_1+\mathrm{A_1}\nu\left(n_s{\bmi q}_1\bcdot{\hat{\bmi z}}- q_sn_1\right)\Big)-\tau_H
L_s n_1, \label{eq34} 
\eeq  
subject to the boundary conditions

\begin{eqnarray}
L_1(x,y,1,\xi,\eta,\nu)&=&0,\,\pi/2\le \theta\le \pi,\, 0\le\phi\le 2\pi,
\label{eq35}\\
L_1(x,y,0,\xi,\eta,\nu)&=&0,\, 0\le \theta\le \pi/2, \, 0\le\phi\le 2\pi.
\label{eq36}
\end{eqnarray} 

The form of Eq. (\ref{eq34}) suggests the following expression for $L_1$:

\beq
L_1=\Psi(z,\xi,\eta,\nu) \exp(\gamma t + i(lx+my)). \label{eq37}   
\eeq

From Eq. (\ref{eq31}) we get

$$
G_1 = \mathcal{G}_1(z)\exp(\gamma t + i(lx+my)),
$$

where

\beq
\mathcal{G}_1(z)=\int_0^{4\pi}\Psi(z,\xi,\eta,\nu)\,\ud\Omega. \label{eq38}
\eeq

Similarly, from Eq. (\ref{eq32}) we get

$$
{\bmi q}_1=\left[q_1^x,q_1^y,q_1^z\right]=\left[P(z),Q(z),S(z)\right] \exp(\gamma t + i(lx+my)), 
$$

where 

\beq
 \left[P(z),Q(z),S(z)\right] =
 \int_0^{4\pi}\big(\xi,\eta,\nu\big)\Psi(z,\xi,\eta,\nu)\,\,\ud\Omega.\label{eq39}  
\eeq

Now from Eqs. (\ref{eq34})--(\ref{eq36}), $\Psi$ satisfies
\beq
\frac{\p \Psi}{\p z} + \frac{i(l\xi+m\eta)+\tau_H n_s}{\nu}\Psi =
\frac{\omega \tau_H}{4\pi\nu}\big(n_s \mathcal{G}_1 + G_s\Theta+\mathrm{A_1}\nu\left(n_s S-
q_s\Theta\right)\big)-\frac{\tau_H}{\nu}L_s\Theta, \label{eq40}
\eeq
subject to the boundary conditions

\begin{eqnarray} 
\Psi(1,\xi,\eta,\nu)=0,\qquad \pi/2\le \theta\le \pi,\quad 0\le\phi\le
2\pi,\label{eq41}\\
\Psi(0,\xi,\eta,\nu)=0,\qquad 0\le \theta\le \pi/2,\quad 0\le\phi\le
2\pi.\label{eq42}    
\end{eqnarray} 

Finally, the linear stability equations become 
\begin{eqnarray}
&& \left(\gamma S_c^{-1}+k^2-\frac{\ud^2}{\ud z^2}\right)\left(\frac{\ud^2}{\ud
z^2}-k^2\right)W=R k^2 \Theta,\label{eq43}\\[6pt]
&& \left(\gamma+k^2-\frac{\ud^2}{\ud z^2}\right)\Theta + V_c \frac{\ud}{\ud
z}\left(M_s\Theta+n_s\frac{\ud M_s}{\ud G_s}\mathcal{G}_1\right)-i\frac{V_c n_s
M_s}{q_s}(lP+mQ) = -\frac{\ud n_s}{\ud z}W,\qquad\quad\label{eq44} 
\end{eqnarray}
subject to the boundary conditions

\begin{eqnarray} 
 W=\frac{\ud W}{\ud z}=\frac{\ud\Theta}{\ud z}-V_c M_s \Theta - V_c n_s
 \frac{\ud M_s}{\ud G_s}
G&=&0, \, \mbox{at}\, z=0, \label{eq45}  \\ [4pt]
W=\frac{\ud^2 W}{\ud z^2}=\frac{\ud\Theta}{\ud z}-V_c M_s \Theta - V_c n_s
\frac{\ud M_s}{\ud G_s}
G&=&0, \, \mbox{at}\, z=1. \label{eq46} 
\end{eqnarray}

Here, $k=\sqrt{l^2+m^2}$ represents the nondimensional wavenumber. Eqs.
(\ref{eq43})--(\ref{eq46}) constitute an eigenvalue problem for $\gamma,$ which is a
function  $\theta_{i}\left(\theta_{0}\right)$, $B,l,m,S_c,V_c,\tau_H,\omega$, $\mathrm{A_1}$ and
$R$ respectively.  If $\mathrm{Re}(\gamma)>0,$ then the equilibrium state is called unstable.

\subsection{Method of solution}

We write $\Psi$ as  $\Psi=\Psi^c+\Psi^d$ to solve Eq. (\ref{eq40}) \cite{mkp20}. Thus,
\begin{eqnarray}
\frac{\p \Psi^c}{\p z} + \frac{i(l\xi+m\eta)+\tau_H n_s}{\nu}\Psi^c &=&
-\frac{\tau_H}{\nu}L_s^c\Theta,\label{eq47}\\[4pt]
\frac{\p \Psi^d}{\p z} + \frac{i(l\xi+m\eta)+\tau_H n_s}{\nu}\Psi^d &=&
\frac{\omega \tau_H}{4\pi\nu}\big(n_s \mathcal{G}_1 + G_s\Theta+\mathrm{A_1}\nu\left(n_s
S-q_s\Theta\right)\big)-\frac{\tau_H}{\nu}L_s^d\Theta,  \qquad
\label{eq48} 
\end{eqnarray}
subject to the boundary conditions
\begin{eqnarray}
\Psi^c(1,\xi,\eta,\nu)=\Psi^d(1,\xi,\eta,\nu)=0,\, \pi/2\le \theta\le
\pi,\, 0\le \phi\le 2\pi,\label{eq49}\\
\Psi^c(0,\xi,\eta,\nu)=\Psi^d(0,\xi,\eta,\nu)=0,\, 0\le \theta\le
\pi/2, \, 0\le \phi\le 2\pi.\label{eq50} 
\end{eqnarray}  
The linear Eq. (\ref{eq47}) is easily solved analytically. Again, we write
$\mathcal{G}_1=\mathcal{G}_1^c+\mathcal{G}_1^d.$ Thus,

$$
\mathcal{G}_1^c=\int_0^{4\pi}\Psi^c\,\ud\Omega=
\tau_H/\cos{\theta_{0}} \left(\int_1^z\Theta(z')\,\ud z' \right)\exp\left(-\tau_H/\cos{\theta_{0}}\int_z^1
n_s(z')\,\ud z'\right),$$ and 
$$  \mathcal{G}_1^d=\int_0^{4\pi}\Psi^d\,\ud\Omega.$$

The integro-differential equation i.e. Eq. (\ref{eq48}) is solved via an iterative method.  $P(z)$ and $Q(z)$ in Eq. (\ref{eq44})
become

\begin{equation}\label{aniso1}
\begin{split}
P(z) & = \int_0^{4\pi}\Psi^d(z,\xi,\eta,\nu)\,\xi\,\ud\Omega,\\
Q(z) & = \int_0^{4\pi}\Psi^d(z,\xi,\eta,\nu)\,\eta\,\ud\Omega. 
\end{split}
\end{equation}

Note that $P(z)$ and $Q(z)$ appear here as a result of scattering only. However, $S(z)$ 
in Eq. (\ref{eq48}) is a function of both the collimated and diffused intensity.

\begin{equation*}
S(z)=\int_0^{4\pi}\Big(\Psi^c(z,\xi,\eta,\nu)
+\Psi^d(z,\xi,\eta,\nu)\Big)\,\nu\,\ud\Omega=-\mathcal{G}_1^{c}+\int_0^{4\pi}\Psi^d(z,\xi,\eta,\nu)\,\nu\,\ud\Omega.
\end{equation*}

We introduce a new variable 

\beq
\widetilde{\Theta}(z)=\int_1^z \Theta(z')\,\ud z'. \label{eq51}
\eeq

Thus, the linear stability equations become  (writing $D=\ud/\ud z$) 
\begin{eqnarray}
&& \left(\gamma S_c^{-1}+k^2-D^2\right)\left(D^2-k^2\right)W=R k^2
D \widetilde{\Theta},\label{eq52}\\[3pt]
&&\Gamma_0(z)+\Gamma_1(z)\widetilde{\Theta} +
\big(\gamma+k^2+\Gamma_2(z)\big)D \widetilde{\Theta} +
V_c M_s\,D^2\widetilde{\Theta}  - D^3 \widetilde{\Theta} =-Dn_s W, \label{eq53} 
\end{eqnarray} 
where 

\begin{eqnarray}
\Gamma_0(z)&=&V_c\,D\left(n_s\frac{\ud M_s}{\ud G_s}\mathcal{G}_1^d\right)-i\frac{V_cn_sM_s}{q_s}(lP+mQ),\label{eq59}\\
\Gamma_1(z)&=& \tau_H/\cos{\theta_{0}} V_c D\left(n_sG_s^c \frac{\ud M_s}{\ud G_s}\right), \label{eq60} \\
\Gamma_2(z)&=& 2\tau_H/\cos{\theta_{0}} V_c n_s G_s^c \frac{\ud M_s}{\ud G_s}+V_c
\frac{\ud M_s}{\ud G_s}DG_s^d.\label{eq61}
\end{eqnarray}  
           
The boundary conditions become
\begin{eqnarray} 
 W=DW=D^2\widetilde{\Theta}-V_c M_s D\widetilde{\Theta} - V_c n_s \frac{\ud
 M_s}{\ud G_s}
 G&=&0, \, \mbox{at}\qquad z=0, \label{eq54}  \\[4pt]
W=D^2W=D^2\widetilde{\Theta}-V_c M_s D\widetilde{\Theta} - V_c n_s \frac{\ud M_s}{\ud
G_s} G&=&0, \, \mbox{at}\qquad z=1, \label{eq55} 
\end{eqnarray}  
and there is an extra boundary condition

\beq
\widetilde{\Theta}(z)=0,\qquad \mbox{at}\quad z=1, \label{eq56}
\eeq

which follows from Eq. (\ref{eq51}).

\section{SOLUTION PROCEDURE}
\label{sec5}

A fourth-order accurate, finite-difference scheme based on Newton–Raphson–Kantorovich (NRK) iterations \cite{cm} is utilized to solve  Eqs. (\ref{eq52}) and (\ref{eq53}) and to calculate the neutral (marginal) stability curves in the $(k, R)$-plane or the growth rate, Re$(\gamma)$, as a function of $R$ for a fixed set of other parameters. The graph of points for which Re$(\gamma) = 0$ is called a marginal (neutral) stability curve. Similarly, if the condition Im$(\gamma) = 0$ satisfies on such a curve, then the bioconvective solution is called stationary (non-oscillatory). Oscillatory solutions exist if Im$(\gamma) \ne 0.$ In addition, overstability occurs if the most unstable mode remains on the oscillatory branch of the neutral curve. When oscillatory solution occurs, a common point (say, $k_b$) exists between the stationary and oscillatory branch. In this instance, the oscillatory branch is the locus of points for which $k \le k_b$. We select here a particular most unstable mode i.e. $\left(k_c, R_c\right)$ of the neutral curve $R^{(n)}(k), \left(n = 1 , 2 , 3 , \cdots \right)$ and in this instance, the wavelength of the initial disturbance is calculated as $\lambda_c=2\pi/k_c$. A bioconvective solution is called mode $n$ if $n$ convection cells can be organized such that one overlies another vertically \cite{mkp20}. Also, the estimated parameters for the proposed problem are same as taken by \cite{mkp20,vh:hv,HPK,gp}.

\section{NUMERICAL RESULTS}
\label{sec6}

We employ the following discrete set of fixed parameters to predict the most unstable mode from an initial equilibrium solution: 
$S_c=20,$  $G_c=1.3,$ $V_c=10,15,20$  $\omega=0.4,$ $B=0.26,0.48,$ $\theta_{i}=0,40,80,$ and $\tau_{H}=0.5, 1,$
at the variation of $\mathrm{A_1}$ such as $\mathrm{A_1}=0,0.4,0.8$. The discrete values of the angle of incidence i.e. $\theta_{i}=0,40,80$ are selected so that the sublayer location at equilibrium state occurs at around mid-height, three-quarter height and top respectively to study the effects of forward scattering on bioconvection.  

The height of the unstable zone (hereafter referred to as HUZ), is defined as the distance of the location of the sublayer at equilibrium state from the bottom of the domain. The difference between the maximum concentration and concentration at the bottom of the suspension is defined as the concentration difference in the unstable zone (hereafter referred to as CDUZ). The maximum concentration is approximately equal to the CDUZ, since in the limit the cell concentration at bottom is zero. A higher value of HUZ or CDUZ implies the suspension is more unstable.

\subsection{Weak forward scattering algal suspension} 
\label{di}

\subsubsection{$V_c=15$} 
\label{}

\section*{Extinction coefficient $\tau_{H}=0.5$} 
\label{ss}

\begin{figure}[!h]
\begin{center}
\includegraphics[scale=0.9]{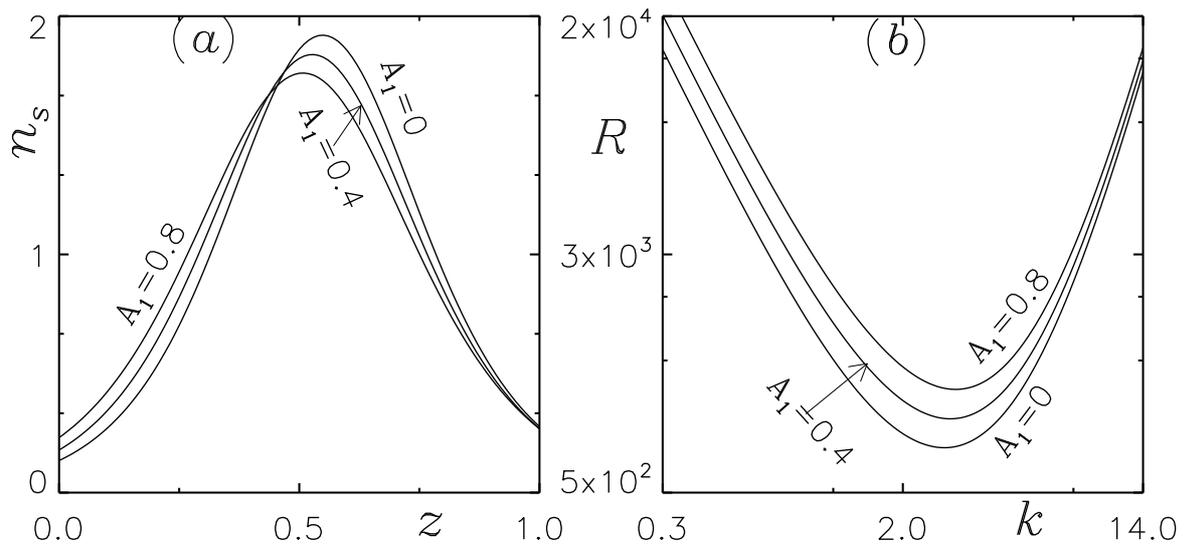}
\end{center}
\caption{(a) Effects of the linearly-anisotropic scattering coefficient $\mathrm{A_1}$ on the base concentration profiles and (b) the corresponding neutral curves. Fixed parameter values are cell swimming speed $V_c=15,$ single scattering albedo $\omega=0.4,$ optical depth $\tau_H=0.5,$ diffuse irradiation  $B=0.26,$ $\theta_{i}=0,$ and critical total intensity $G_{c}=1.3$.}
\label{fig6}
\end{figure}

Figure \ref{fig6} depicts the effects of the forward scattering coefficient $\mathrm{A_1}$ on bioconvective instability, when the parameters $V_c=15,$ $\tau_{H}=0.5,$ $\omega=0.4,$ $\theta_{i}=0,$ and $B=0.26$ are kept fixed. Figures \ref{fig6}(a) and \ref{fig6}(b) show the  equilibrium solutions  at sublayer location and the corresponding unstable modes. For $\mathrm{A_1}=0,$ the sublayer at equilibrium state forms at around mid-height of the domain and the location of sublayer at equilibrium state shifts towards the bottom as $\mathrm{A_1}$ is increased to $0.4$ and $0.8$ respectively. Hence, both the HUZ and CDUZ decrease as $\mathrm{A_1}$ is increased from $0$ to $0.8.$ Thus, the critical Rayleigh number $R_c$ increases as $\mathrm{A_1}$ is increased from $0$ to $0.8$ making the suspension more stable [see Fig.~\ref{fig6}].

\begin{figure}[!h]
\begin{center}
\includegraphics[scale=0.9]{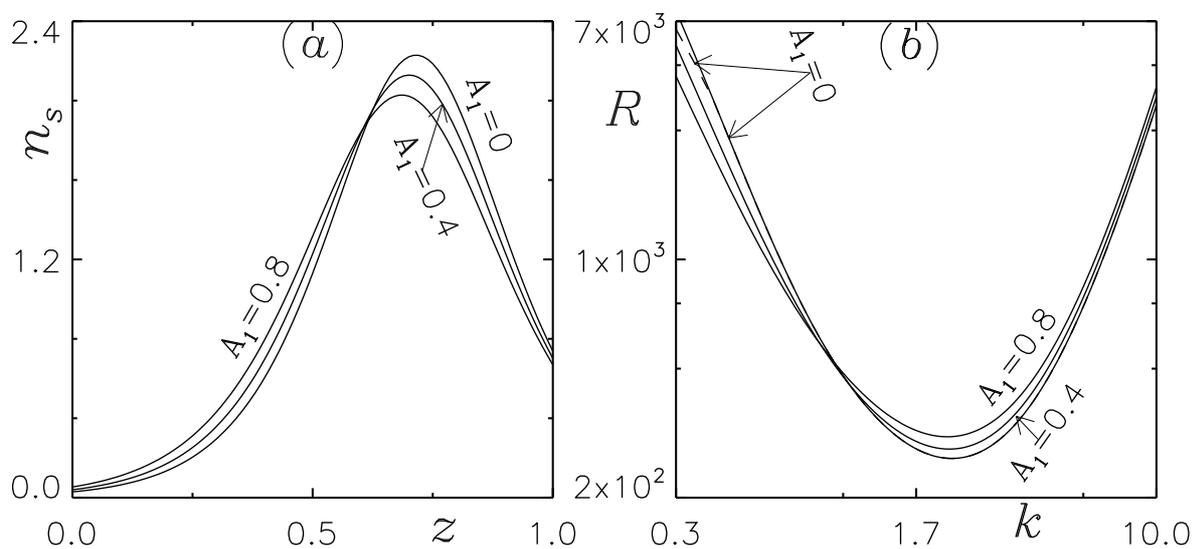}
\end{center}
\caption{(a) Effects of the linearly-anisotropic scattering coefficient $\mathrm{A_1}$ on the base concentration profiles and (b) the corresponding neutral curves. Fixed parameter values are cell swimming speed $V_c=15,$ single scattering albedo $\omega=0.4,$ optical depth $\tau_H=0.5,$  diffuse irradiation $B=0.26,$ $\theta_{i}=40,$ and critical total intensity $G_{c}=1.3$.}
\label{fig7}
\end{figure}

Now consider the case when $V_c=15,$ $\tau_{H}=0.5,$ $\omega=0.4,$ $\theta_{i}=40,$ and $B=0.26.$ In this instance, when $\mathrm{A_1}=0,$ the sublayer at equilibrium state forms at around three-quarter height of the domain and the location of sublayer at equilibrium state shifts towards the bottom as $\mathrm{A_1}$ is increased to $0.4$ and $0.8$ respectively. Hence, both the HUZ and CDUZ decrease as $\mathrm{A_1}$ is increased from $0$ to $0.8.$ Thus, the critical Rayleigh number $R_c$ increases as $\mathrm{A_1}$ is increased from $0$ to $0.8$ making the suspension more stable [see Fig.~\ref{fig7}].

\begin{figure}[!h]
\begin{center}
\includegraphics[scale=0.9]{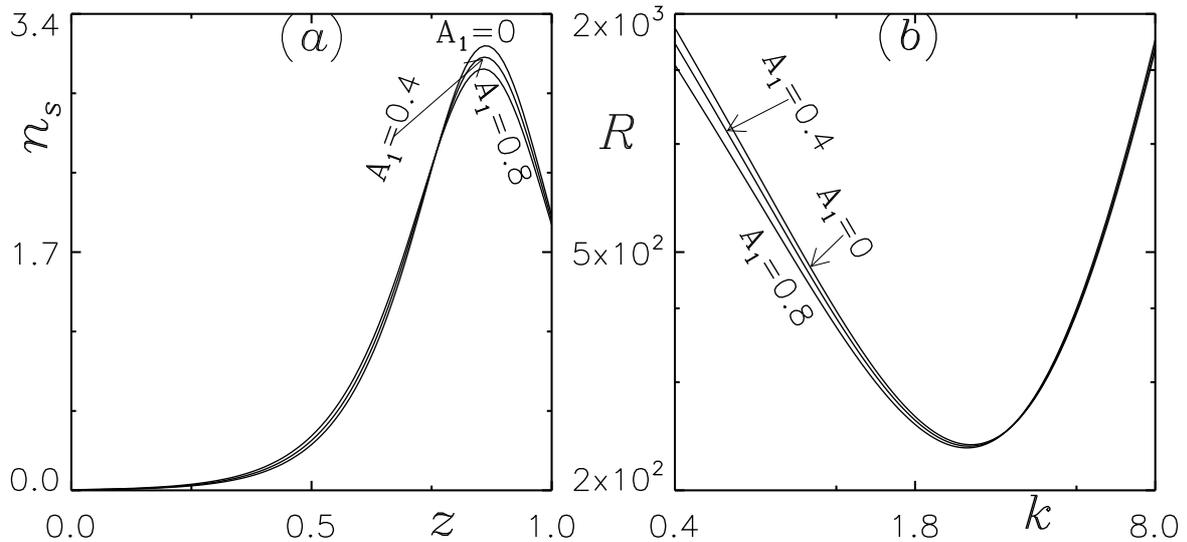}
\end{center}
\caption{(a) Effects of the linearly-anisotropic scattering coefficient $\mathrm{A_1}$ on the base concentration profiles and (b) the corresponding neutral curves. Fixed parameter values are cell swimming speed $V_c=15,$ single scattering albedo $\omega=0.4,$ optical depth $\tau_H=0.5,$  diffuse irradiation $B=0.26,$ $\theta_{i}=80,$ and critical total intensity $G_{c}=1.3$.}
\label{fig8}
\end{figure}

Next, consider the case when $V_c=15,$ $\tau_{H}=0.5,$ $\omega=0.4,$ $\theta_{i}=80,$ and $B=0.26.$ In this instance, when $\mathrm{A_1}=0,$ the sublayer at equilibrium state forms at around top of the domain and the location of sublayer at equilibrium state is same as $\mathrm{A_1}$ is increased to $0.4$ and $0.8$ respectively. Here the HUZ is almost same for $\mathrm{A_1}=0,0.4,0.8.$ But CDUZ decreases as $\mathrm{A_1}$ is increased from $0$ to $0.8.$ Here the latter effect is more dominant than the former and thus, the critical Rayleigh number $R_c$ decreases as $\mathrm{A_1}$ is increased from $0$ to $0.8$ making the suspension more stable [see Fig.~\ref{fig8}]. The most unstable mode from an initial equilibrium solution is stationary for $\theta_{i}=0,40,80$ when $\tau_{H}=0.5.$

\begin{figure}[!h]
\begin{center}
\includegraphics[scale=0.9]{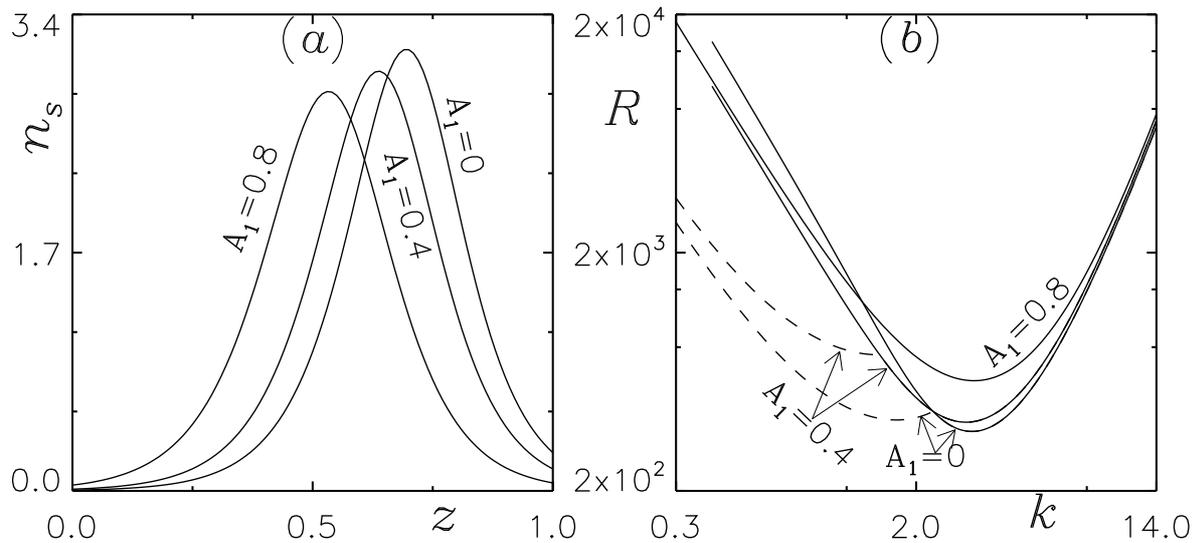}
\end{center}
\caption{(a) Effects of the linearly-anisotropic scattering coefficient $\mathrm{A_1}$ on the base concentration profiles and (b) the corresponding neutral curves. Fixed parameter values are cell swimming speed $V_c=15,$ single scattering albedo $\omega=0.4,$ optical depth $\tau_H=1,$  diffuse irradiation $B=0.48,$ $\theta_{i}=0,$ and critical total intensity $G_{c}=1.3$.}
\label{fig9}
\end{figure}

Now consider the case when $V_c=15,$ $\tau_{H}=1,$ $\omega=0.4,$ $\theta_{i}=0$ and $B=0.48,$ and Fig. \ref{fig9} shows the effects of the forward scattering coefficient $\mathrm{A_1}$ on bioconvective instability. In this instance, when $\mathrm{A_1}=0,$ the sublayer at equilibrium state forms at around three-quarter height of the domain and the location of sublayer at equilibrium state shifts towards the bottom as $\mathrm{A_1}$ is increased to $0.4$ and $0.8$ respectively. Hence, both the HUZ and CDUZ decrease as $\mathrm{A_1}$ is increased from $0$ to $0.8.$ Thus, the critical Rayleigh number $R_c$ increases as $\mathrm{A_1}$ is increased from $0$ to $0.8$ making the suspension more stable [see Fig.~\ref{fig9}]. Also, the most unstable mode from an initial equilibrium solution is stationary for $\mathrm{A_1}=0,0.4.$ However, an oscillatory branch bifurcates from the stationary branch in both the cases. When $\mathrm{A_1}$ is increased to $0.8,$ the oscillarory branch disappears [see Fig.~\ref{fig9}].

\section*{Extinction coefficient $\tau_{H}=1.0$} 
\label{os}

\begin{figure}[!h]
\begin{center}
\includegraphics[scale=0.9]{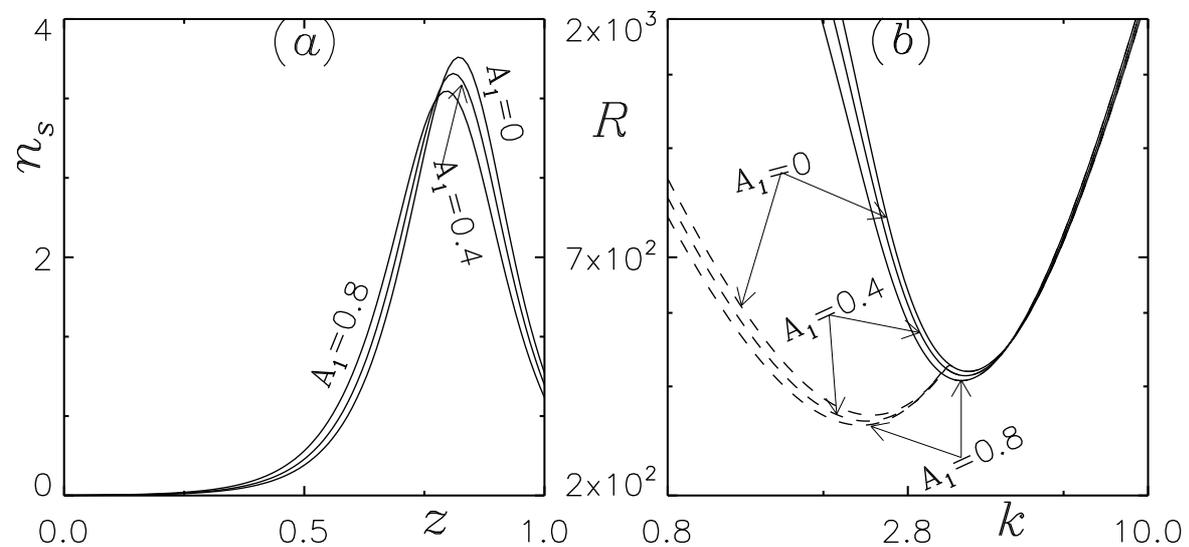}
\end{center}
\caption{(a) Effects of the linearly-anisotropic scattering coefficient $\mathrm{A_1}$ on the base concentration profiles and (b) the corresponding neutral curves. Fixed parameter values are cell swimming speed $V_c=15,$ single scattering albedo $\omega=0.4,$ optical depth $\tau_H=1,$  diffuse irradiation $B=0.48,$ $\theta_{i}=40,$ and critical total intensity $G_{c}=1.3$.}
\label{fig10}
\end{figure}

Now consider the case when $V_c=15,$ $\tau_{H}=1,$ $\omega=0.4,$ $\theta_{i}=40$ and $B=0.48,$ and Fig. \ref{fig10} shows the effects of the forward scattering coefficient $\mathrm{A_1}$ on bioconvective instability. In this instance, when $\mathrm{A_1}=0,$ the sublayer at equilibrium state forms at around $z \approx 0.8$ of the domain and the location of sublayer at equilibrium state shifts towards the bottom as $\mathrm{A_1}$ is increased to $0.4$ and $0.8$ respectively. The most unstable mode from an initial equilibrium solution is ovrstable for $\mathrm{A_1}=0,0.4,0.8$ [see Fig.~\ref{fig10}].

\begin{figure}[!h]
\begin{center}
\includegraphics[scale=0.9]{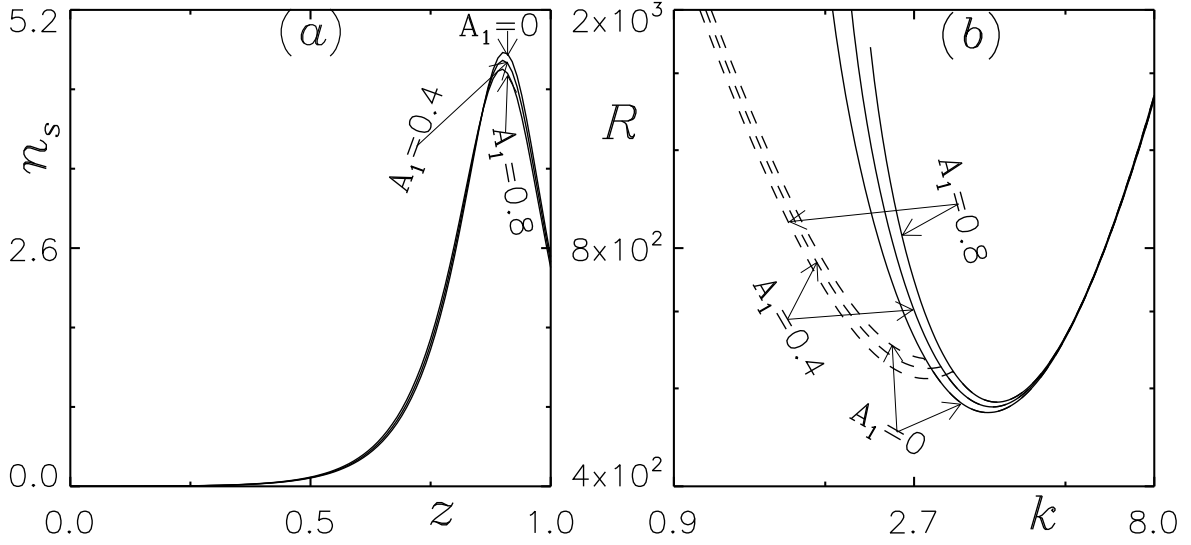}
\end{center}
\caption{(a) Effects of the linearly-anisotropic scattering coefficient $\mathrm{A_1}$ on the base concentration profiles and (b) the corresponding neutral curves. Fixed parameter values are cell swimming speed $V_c=15,$ single scattering albedo $\omega=0.4,$ optical depth $\tau_H=1,$  diffuse irradiation $B=0.48,$ $\theta_{i}=80,$ and critical total intensity $G_{c}=1.3$.}
\label{fig11}
\end{figure}

\begin{table}[!h]
\begin{center}
\caption{Bioconvective solutions with the variation $\mathrm{A},$ by keeping other governing parameters fixed. A  result with double dagger symbol indicates that a smaller minimum occurs on an oscillatory branch and a starred result indicates that $R^{(1)}(k)$ branch of the neutral curve is oscillatory. \label{nr1}}
\begin{tabular}{c c c c c c c c c}
\hline \hline			
$V_c$ \qquad& $\tau_{H}$ \qquad& $\omega$ \qquad& $B$ \qquad& $\theta_i$ \qquad& $A_1$ & $\lambda_c$ \qquad& $R_c$ \qquad& $Im(\gamma)$ \\[0.25pt] \hline
15 \qquad& 0.5 \qquad& 0.4 \qquad& 0.26 \qquad& 0 \qquad& 0 \qquad& 2.20 \qquad& 709.68 \qquad& 0 \\[0.25pt]
15 \qquad& 0.5 \qquad& 0.4 \qquad& 0.26 \qquad& 0 \qquad& 0.4 \qquad& 2.11 \qquad& 889.06 \qquad& 0 \\[0.25pt]
15 \qquad& 0.5 \qquad& 0.4 \qquad& 0.26 \qquad& 0 \qquad& 0.8 \qquad& 2 \qquad& 1118.07 \qquad& 0 \\[0.25pt]
15 \qquad& 0.5 \qquad& 0.4 \qquad& 0.26 \qquad& 40 \qquad& $0^{\star}$ \qquad& 2.81 \qquad& 266.34 \qquad& 0 \\[0.25pt]
15 \qquad& 0.5 \qquad& 0.4 \qquad& 0.26 \qquad& 40 \qquad& 0.4 \qquad& 2.88 \qquad& 285.21 \qquad& 0 \\[0.25pt]
15 \qquad& 0.5 \qquad& 0.4 \qquad& 0.26 \qquad& 40 \qquad& 0.8 \qquad& 2.88 \qquad& 312.01 \qquad& 0 \\[0.25pt]
15 \qquad& 0.5 \qquad& 0.4 \qquad& 0.26 \qquad& 80 \qquad& 0 \qquad& 2.49 \qquad& 242.57 \qquad& 0 \\[0.25pt]
15 \qquad& 0.5 \qquad& 0.4 \qquad& 0.26 \qquad& 80 \qquad& 0.4 \qquad& 2.49 \qquad& 241 \qquad& 0 \\[0.25pt]
15 \qquad& 0.5 \qquad& 0.4 \qquad& 0.26 \qquad& 80 \qquad& 0.8 \qquad& 2.55 \qquad& 239.32 \qquad& 0 \\	[0.25pt]		
15 \qquad& 1 \qquad& 0.4 \qquad& 0.48 \qquad& 0 \qquad& $0^{\star}$ \qquad& 1.96 \qquad& 356.29 \qquad& 0 \\[0.25pt]
15 \qquad& 1 \qquad& 0.4 \qquad& 0.48 \qquad& 0 \qquad& $0.4^{\star}$ \qquad& 2.07 \qquad& 389.90 \qquad& 0 \\[0.25pt]
15 \qquad& 1 \qquad& 0.4 \qquad& 0.48 \qquad& 0 \qquad& 0.8 \qquad& 1.93 \qquad& 583.55 \qquad& 0 \\[0.25pt]
15 \qquad& 1 \qquad& 0.4 \qquad& 0.48 \qquad& 40 \qquad& $0^{\ddag}$ \qquad& 2.67 \qquad& 354.98 \qquad& 12.57 \\[0.25pt]
15 \qquad& 1 \qquad& 0.4 \qquad& 0.48 \qquad& 40 \qquad& $0.4^{\ddag}$ \qquad& 2.74 \qquad& 345.88 \qquad& 12.14 \\[0.25pt]
15 \qquad& 1 \qquad& 0.4 \qquad& 0.48 \qquad& 40 \qquad& $0.8^{\ddag}$ \qquad& 2.88 \qquad& 339.71 \qquad& 11.58 \\[0.25pt]
15 \qquad& 1 \qquad& 0.4 \qquad& 0.48 \qquad& 80 \qquad& $0^{\star}$ \qquad& 1.67 \qquad& 490.81 \qquad& 0 \\[0.25pt]
15 \qquad& 1 \qquad& 0.4 \qquad& 0.48 \qquad& 80 \qquad& $0.4^{\star}$ \qquad& 1.64 \qquad& 498.58 \qquad& 0 \\[0.25pt]
15 \qquad& 1 \qquad& 0.4 \qquad& 0.48 \qquad& 80 \qquad& $0.8^{\star}$ \qquad& 1.60 \qquad& 505.28 \qquad& 0 \\[0.25pt]
\hline \hline
\end{tabular}
\end{center}
\end{table}

Now consider the case when $V_c=15,$ $\tau_{H}=1,$ $\omega=0.4,$ $\theta_{i}=80$ and $B=0.48,$ and Fig. \ref{fig11} shows the effects of the forward scattering coefficient $\mathrm{A_1}$ on bioconvective instability. In this instance, when $\mathrm{A_1}=0,$ the sublayer at equilibrium state forms at around $z \approx 0.9$ of the domain and the location of sublayer at equilibrium state shifts towards the bottom as $\mathrm{A_1}$ is increased to $0.4$ and $0.8$ respectively. The most unstable mode from an initial equilibrium solution is stationary for $\mathrm{A_1}=0,0.4,0.8$ However, an oscillatory branch bifurcates from the stationary branch in the three cases [see Fig.~\ref{fig11}].

The numerical results for the critical Rayleigh number $(R_c)$ and the wavelength $(\lambda_c)$ of this section are summarized in Table~\ref{nr1}.

\subsubsection{$V_c=10$ \textnormal{and}  $20$} 
\label{}

\begin{table}[!h]
\begin{center}
\caption{Bioconvective solutions with the variation  $\mathrm{A},$ by keeping other governing parameters fixed. A  result with double dagger symbol indicates that a smaller minimum occurs on an oscillatory branch and a starred result indicates that $R^{(1)}(k)$ branch of the neutral curve is oscillatory.  \label{nr2}}
\begin{tabular}{c c c c c c c c c}
\hline \hline
$V_c$ & $\tau_{H}$ & $\omega$ & $B$ & $\theta_i$ & $A_1$ & $\lambda_c$ & $R_c$ & $Im(\gamma)$ \\
[0.25pt] \hline
10 \qquad& 0.5 \qquad& 0.4 \qquad& 0.26 \qquad& 0 \qquad& 0 \qquad& 2.11 \qquad& 1385.54 \qquad& 0 \\[0.25pt]
10 \qquad& 0.5 \qquad& 0.4 \qquad& 0.26 \qquad& 0 \qquad& 0.4 \qquad& 2.07 \qquad& 1635.99 \qquad& 0 \\[0.25pt]
10 \qquad& 0.5 \qquad& 0.4 \qquad& 0.26 \qquad& 0 \qquad& 0.8 \qquad& 2 \qquad& 1949.36 \qquad& 0 \\[0.25pt]
10 \qquad& 0.5 \qquad& 0.4 \qquad& 0.26 \qquad& 40 \qquad& 0 \qquad& 3.21 \qquad& 497.16 \qquad& 0 \\[0.25pt]
10 \qquad& 0.5 \qquad& 0.4 \qquad& 0.26 \qquad& 40 \qquad& 0.4 \qquad& 3.12 \qquad& 565.18 \qquad& 0 \\[0.25pt]
10 \qquad& 0.5 \qquad& 0.4 \qquad& 0.26 \qquad& 40 \qquad& 0.8 \qquad& 3.04 \qquad& 649.09 \qquad& 0 \\[0.25pt]
10 \qquad& 0.5 \qquad& 0.4 \qquad& 0.26 \qquad& 80 \qquad& 0 \qquad& 3.75 \qquad& 200.88 \qquad& 0 \\[0.25pt]
10 \qquad& 0.5 \qquad& 0.4 \qquad& 0.26 \qquad& 80 \qquad& 0.4 \qquad& 3.88 \qquad& 202.47 \qquad& 0 \\[0.25pt]
10 \qquad& 0.5 \qquad& 0.4 \qquad& 0.26 \qquad& 80 \qquad& 0.8 \qquad& 4.02 \qquad& 204.29 \qquad& 0 \\[0.25pt]
10 \qquad& 1 \qquad& 0.4 \qquad& 0.48 \qquad& 0 \qquad& 0 \qquad& 2.24 \qquad& 523.41 \qquad& 0 \\[0.25pt]
10 \qquad& 1 \qquad& 0.4 \qquad& 0.48 \qquad& 0 \qquad& 0.4 \qquad& 2.11 \qquad& 663.57 \qquad& 0 \\[0.25pt]
10 \qquad& 1 \qquad& 0.4 \qquad& 0.48 \qquad& 0 \qquad& 0.8 \qquad& 2 \qquad& 874.99 \qquad& 0 \\[0.25pt]
10 \qquad& 1 \qquad& 0.4 \qquad& 0.48 \qquad& 40 \qquad& $0^{\star}$ \qquad& 2.29 \qquad& 316.56 \qquad& 0 \\[0.25pt]
10 \qquad& 1 \qquad& 0.4 \qquad& 0.48 \qquad& 40 \qquad& $0.4^{\star}$ \qquad& 2.29 \qquad& 328.63 \qquad& 0 \\[0.25pt]
10 \qquad& 1 \qquad& 0.4 \qquad& 0.48 \qquad& 40 \qquad& $0.8^{\star}$ \qquad& 2.34 \qquad& 346.95 \qquad& 0 \\[0.25pt]
10 \qquad& 1 \qquad& 0.4 \qquad& 0.48 \qquad& 80 \qquad& 0 \qquad& 2.24 \qquad& 297.52 \qquad& 0 \\[0.25pt]
10 \qquad& 1 \qquad& 0.4 \qquad& 0.48 \qquad& 80 \qquad& 0.4 \qquad& 2.24 \qquad& 299.77 \qquad& 0 \\[0.25pt]
10 \qquad& 1 \qquad& 0.4 \qquad& 0.48 \qquad& 80 \qquad& 0.8 \qquad& 2.24 \qquad& 301.90 \qquad& 0 \\[0.25pt]
20 \qquad& 0.5 \qquad& 0.4 \qquad& 0.26 \qquad& 0 \qquad& $0^{\star}$ \qquad& 2.29 \qquad& 413.01 \qquad& 0 \\[0.25pt]
20 \qquad& 0.5 \qquad& 0.4 \qquad& 0.26 \qquad& 0 \qquad& 0.4 \qquad& 2.15 \qquad& 545.37 \qquad& 0 \\[0.25pt]
20 \qquad& 0.5 \qquad& 0.4 \qquad& 0.26 \qquad& 0 \qquad& 0.8 \qquad& 2 \qquad& 745.32 \qquad& 0 \\[0.25pt]
20 \qquad& 0.5 \qquad& 0.4 \qquad& 0.26 \qquad& 40 \qquad& $0^{\star}$ \qquad& 2.2 \qquad& 274.16 \qquad& 0 \\[0.25pt]
20 \qquad& 0.5 \qquad& 0.4 \qquad& 0.26 \qquad& 40 \qquad& $0.4^{\star}$ \qquad& 2.24 \qquad& 275.25 \qquad& 0 \\[0.25pt]
20 \qquad& 0.5 \qquad& 0.4 \qquad& 0.26 \qquad& 40 \qquad& 0.8 \qquad& 2.34 \qquad& 279.11 \qquad& 0 \\[0.25pt]
20 \qquad& 0.5 \qquad& 0.4 \qquad& 0.26 \qquad& 80 \qquad& $0^{\star}$ \qquad& 1.9 \qquad& 355.66 \qquad& 0 \\[0.25pt]
20 \qquad& 0.5 \qquad& 0.4 \qquad& 0.26 \qquad& 80 \qquad& $0.4^{\star}$ \qquad& 1.93 \qquad& 352.28 \qquad& 0 \\[0.25pt]
20 \qquad& 0.5 \qquad& 0.4 \qquad& 0.26 \qquad& 80 \qquad& $0.8^{\star}$ \qquad& 1.93 \qquad& 348.54 \qquad& 0 \\[0.25pt]
20 \qquad& 1 \qquad& 0.4 \qquad& 0.48 \qquad& 0 \qquad& $0^\ddag$ \qquad& 3.31 \qquad& 266.52 \qquad& 14.62 \\[0.25pt]
20 \qquad& 1 \qquad& 0.4 \qquad& 0.48 \qquad& 0 \qquad& 0.4 \qquad& 2.11 \qquad& 663.57 \qquad& 0 \\[0.25pt]
20 \qquad& 1 \qquad& 0.4 \qquad& 0.48 \qquad& 0 \qquad& 0.8 \qquad& 1.9 \qquad& 460.76 \qquad& 0 \\[0.25pt]
20 \qquad& 1 \qquad& 0.4 \qquad& 0.48 \qquad& 40 \qquad& $0^\ddag$ \qquad& 2.44 \qquad& 400.79 \qquad& 23.95 \\[0.25pt]
20 \qquad& 1 \qquad& 0.4 \qquad& 0.48 \qquad& 40 \qquad& $0.4^\ddag$ \qquad& 2.49 \qquad& 375.50 \qquad& 23.04 \\[0.25pt]
20 \qquad& 1 \qquad& 0.4 \qquad& 0.48 \qquad& 40 \qquad& $0.8^\ddag$ \qquad& 2.55 \qquad& 351.02 \qquad& 21.90 \\[0.25pt]
20 \qquad& 1 \qquad& 0.4 \qquad& 0.48 \qquad& 80 \qquad& $0^\ddag$ \qquad& 1.96 \qquad& 793.05 \qquad& 17.46 \\[0.25pt]
20 \qquad& 1 \qquad& 0.4 \qquad& 0.48 \qquad& 80 \qquad& $0.4^\ddag$ \qquad& 1.96 \qquad& 765.71 \qquad& 19.99 \\[0.25pt]
20 \qquad& 1 \qquad& 0.4 \qquad& 0.48 \qquad& 80 \qquad& $0.8^\ddag$ \qquad& 2 \qquad& 737.49 \qquad& 22.18 \\[0.25pt]
\hline \hline
\end{tabular}
\end{center}
\end{table}

The prediction about the most unstable mode from an equilibrium solution at the variation of the forward scattering coeficient at bioconvective instability for $V_c=10$ and $20$ are qualitatively similar to those of $V_c=15$. The numerical results bioconvective instability of this section are addressed in Table~\ref{nr2}.

\section{CONCLUSIONS}
\label{sec7}

The proposed phototaxis model investigates the effect of anisotropic (forward) scattering by algae on a most unstable mode from an initial equilibrium solution at bioconvective instability. Accounting the phototaxis model to be more realistic and reliable for industrial exploitation, the illuminating source of the algal suspension is composed of both diffuse and oblique (i.e. not vertical) collimated flux. Furthermore, it is assumed that the scattering is azimuthally symmetric as we neglect bottom-heaviness in algae. Thus, the anisotropic phase function is linear.

The conclusion regarding the  obtained numerical results is drawn based on discrete parameter range and is outlined as follows. In the equilibrium state, the value of the total intensity in the top (resp. bottom) half of the suspension decreases (resp. increases)  as the forward scattering coefficient $\mathrm{A_1}$  is increased in a uniform suspension. Moreover, the location of the critical value of total intensity shifts towards the bottom (resp. top ) in the  bottom (resp. top) half of the uniform suspension as the forward scattering coefficient $\mathrm{A_1}$  is increased. At basic equilibrium state, the values of HUZ and CDUZ decrease as the forward scattering coefficient $\mathrm{A_1}$  is increased. Thus, location of sublayer at equilibrium state shifts toward the bottom of the domain as the forward scattering coefficient $\mathrm{A_1}$  is increased. 

In case of stationary bioconvective solutions, the critical Rayleigh number increases at instability as the forward scattering coefficient $\mathrm{A_1}$  is increased as the values of HUZ and CDUZ decrease in the equilibrium state. This often happens when the sublayer at equilibrium state occurs at around mid-height and three-quarter height of the domain based on the value of angle of incidence $\theta_{i}.$ In contrast, the critical Rayleigh number decreases at instability as the forward scattering coefficient $\mathrm{A_1}$  is increased when the sublayer at equilibrium state forms at around the top.

Next, we analyse about the oscillatory bioconvective solutions. In this instance, the oscillatory branches which appear at $\mathrm{A_1}=0,0.4$ suddenly disappear as $\mathrm{A_1}$ is increased to $0.8.$ This is due to the sublayer at equilibrium state shifts toward the bottom as $\mathrm{A_1}$ is increased to $0.8.$ Similarly, when the sublayer at equilibrium state forms at around three-quarter height of the domain, the most unstable mode remains in the oscillatory branch making the bioconvective solution overstable for all values of $\mathrm{A_1}$ i.e. when $\mathrm{A_1}=0,0.4,0.8.$ Again, as the sublayer at equilibrium state forms at around top of the domain, an oscillatory branch bifurcates from the stationary branch. However, the most unstable mode still remains in the stationary branch for all values of $\mathrm{A_1}$ i.e. when $\mathrm{A_1}=0,0.4,0.8.$

Since algae are mainly predominantly photo-gyro-gravitactic in a natural environment despite of being purely phototactic \cite{williams_11,hh_1997,khh:hhk}. Hence, it is not possible to find a quantitative analysis of similar experiments to compare with our theoretical developments via the proposed model. Panda \textit{et al.} \cite{prmm16} used an vertical (i.e. not oblique) collimated flux in addition to a diffuse flux as an illuminating source to the algal suspension in their study. Finally, the sunlight strikes an algal suspension mostly at off-normal angles via an oblique collimated flux and thus, the proposed model is more realistic in nature than Panda \textit{et al.} \cite{prmm16}.

\section*{Acknowledgements}

The corresponding author gratefully acknowledges the Science $\&$ Engineering Research Board (SERB), a statutory body of the Department of Science $\&$ Technology (DST), Government of India for the financial support via the Core Research Grant (Grant No. CRG/2018/000105).
Also, the second author gratefully acknowledges the Ministry of Education (Government of India) for the financial support via
GATE fellowship (Registration No. MA19S43047204).

\bibliographystyle{plain}

\end{document}